\documentclass[a4paper,11pt]{amsart} 
\usepackage{amscd}             
\usepackage{amssymb}
\usepackage{amsthm}
\usepackage{epsf} 
\usepackage[T1]{fontenc}            
\newtheorem{thm}{Theorem}[section]   
\newtheorem{cor}[thm]{Corollary}
\newtheorem{lemma}[thm]{Lemma}
\newtheorem{prop}[thm]{Proposition}

\newtheorem{rem}[thm]{Remark}

\def\c1{\operatorname{c_1}}

\def\ZZ{{\mathbf Z}}

\def\PP{{\mathbf P}}

\def\O{{\mathcal O}}

\def\iso{\simeq}                 
\def\eqv{\equiv}

\def\+{\oplus}                   
\def\*{\otimes}                  
\def\hpil{\longrightarrow}       
\def\khpil{\rightarrow}

\def\Pic{\operatorname{Pic}}

\def\mod{\operatorname{mod}}
\def\disc{\operatorname{disc}}

\def\k{\Lambda}

\def\hs{\hspace{.05in}}

\begin{document}

\title{Smooth Curves on Projective $K3$ Surfaces}
\author{Andreas Leopold Knutsen}  
\address{Dept. of Mathematics\\ 
  University of Bergen\\ Johs. Brunsgt 12\\ 5008 Bergen \\ Norway}

\begin{abstract}
  In this paper we give for all $n \geq 2$, $d>0$, $g \geq 0$ necessary
  and sufficient conditions for the existence of a pair $(X,C)$, where
  $X$ is a $K3$ surface of degree $2n$ in $\PP^{n+1}$ and $C$ is a
  smooth (reduced and irreducible) curve of degree $d$ and genus $g$ on $X$. 
  The surfaces constructed have Picard group of minimal rank possible 
  (being either 
  $1$ or  $2$), and in each case we specify a set of generators. For 
  $n \geq 4$ we also  determine when $X$ can be chosen to be an intersection 
  of quadrics (in all other cases $X$ has to be an intersection of both 
  quadrics and cubics). Finally, we give necessary and sufficient 
  conditions for $\O_C (k)$ to be non-special, for any integer $k \geq
  1$.
\end{abstract}

\maketitle

\section{Introduction}
\label{intro} 

In recent years the interest for $K3$ surfaces and Calabi-Yau threefolds has 
increased because of their importance in theoretical physics and string 
theory in particular. The study of curves on $K3$ surfaces is interesting not 
only in its own right, but also because one can use $K3$ surfaces containing 
particular curves to constuct $K3$ fibered Calabi-Yau threefolds containing 
the same curves as rigid curves (see \cite{katz}, \cite{oguiso}, \cite{EJS}, 
\cite{kley} and \cite{kley2}).

The problem of determining the possible pairs $(d,g)$
of degree  $d$ and 
genus $g$ of curves contained in certain ambient varieties is rather
fascinating. A fundamental result of L.~Gruson and C.~Peskine in \cite{grus} 
determines all such pairs for which
there exists a smooth irreducible nondegenerate curve of degree  $d$ and
genus $g$ in $\PP^3$. To solve the problem, the authors need curves on
some rational quartic surface with a double line.

S.~Mori proved in \cite{mori} that essentially the same degrees
and genera as those found by Gruson and Peskine for curves on rational quartic
surfaces, can be found on smooth quartic surfaces as well.

K.~Oguiso \cite{oguiso} showed in 1994 that for all $n \geq 2$ and $d>0$
there exists a $K3$ surface of degree $2n$ containing a smooth rational
curve of degree $d$.

The main aim of this paper is to prove the following general result:

\begin{thm}   \label{comp.char}
  Let $n \geq 2$, $d>0$, $g \geq 0$ be integers. Then there exists a
  $K3$ surface \footnote{By a $K3$ surface is meant a smooth $K3$
  surface.}  $X$ of degree $2n$ in $\PP^{n+1}$ containing a smooth 
  curve $C$  of degree $d$ and genus $g$ if and only if 
    \begin{itemize}
     \item[(i)] $g=d^2/4n + 1$ and there exist integers $k,m \geq 1$
       and $(k,m) \not = (2,1)$ such that $n = k^2 m$ and $2n$ divides
       $kd$,   
     \item[(ii)] $d^2/4n < g < d^2/4n +1$ except in the following cases
       \begin{itemize}
        \item[(a)] $d \eqv \pm 1, \pm 2 \hs (\mod 2n)$,
        \item[(b)] $d^2-4n(g-1) = 1$ and $d \eqv n \pm 1 \hs (\mod 2n)$,
        \item[(c)] $d^2-4n(g-1) = n$ and $d \eqv n \hs (\mod 2n)$,
        \item[(d)] $d^2-4n(g-1) = 1$ and $d-1$ or $d+1$ divides $2n$,
        \end{itemize}
\item[(iii)] $g = d^2/4n$ and $d$ is not divisible by $2n$,     
\item[(iv)] $g < d^2/4n$ and $(d,g) \not =(2n+1,n+1)$.
    \end{itemize}
  
  Furthermore, in case (i) $X$ can be chosen such that $\Pic X = \ZZ
  \frac{2n}{dk} C = \ZZ \frac{1}{k}H$ and in cases (ii)-(iv) such that $\Pic
  X = \ZZ H \+ \ZZ C$, where $H$ is the hyperplane section of $X$.

  If $n \geq 4$, $X$ can be chosen to be scheme-theoretically an
  intersection of quadrics in cases (i), (iii) and (iv), and also in case 
  (ii), except when $d^2-4n(g-1)=1$ and $3d \eqv  \pm 3 \hs (\mod 2n)$ or 
  $d^2-4n(g-1)=9$ and $d \eqv  \pm 3 \hs (\mod 2n)$, in which case $X$
  has to be an intersection of both quadrics and cubics.
\end{thm}

\begin{rem} \label{mainrem}
  If one allows $X$ to be a birational projective model of a $K3$
  surface (which automatically yields with at worst rational double points as
  singularities), then the result above remains the same, except that
  the case (ii)-(c) occurs.
\end{rem}

The most general results concerning construction of $K3$ fibered Calabi-Yau 
threefolds are due to H.~P.~Kley \cite{kley2}, who constructs rigid curves of 
bounded genera on complete intersection Calabi-Yau threefolds ($CICY$s). The 
approach of Kley requires that the smooth curve $C$ on the $K3$ surface $X$ 
used to construct the $CICY$ is linearly independent of the hyperplane 
section $H$ of $X$ and also 
that $h^1 (C',\O_{C'} (k))=0$ for all $C'\in |C|$ for $k=1$ or $2$
(depending on the different types of $CICY$s). Motivated by this, we
also prove the 
following result, which is an improvement of the results in \cite{kley2} and 
gives the corresponding existence of more rigid curves in $CICY$s than is 
shown in \cite{kley2}.

\begin{prop}   \label{extra} 
Let $k \geq 1$ be an integer. We can find $X$ and $C$ as in Theorem
\ref{comp.char} such that $h^1 (C', \O_{C'} (k))=0$ for all $C'\in
|C|$ if and only if
\[ d \leq 2nk \hs \mbox{ or } dk > nk^2 +g.  \]
\end{prop}

So far one has only used the $K3$ surfaces that are complete intersections (more specifically the smooth complete intersections of type $(4)$ in  $\PP^3$,
$(2,3)$ in  $\PP^4$ and $(2,2,2)$ in  $\PP^5$, see Section \ref{app})
to construct $CICY$s containing rigid curves. S. Mukai showed in
\cite{mukai} that general $K3$ surfaces of degrees $10$, $12$, $14$,
$16$ and $18$ are complete intersections in homogeneous spaces. For
the triples $(n,d,g)$ in Theorem \ref{comp.char} corresponding to such
general surfaces, one can then construct $K3$ fibered Calabi-Yau
threefolds that are complete intersections in homogeneous spaces
containing rigid curves. This is the topic in \cite{kn2}.

We work over the field of
complex numbers. The results will hold for any
algebraically closed field of characteristic zero by the Lefschetz principle. 

It is a pleasure to thank Professor Trygve Johnsen at the University of 
Bergen. I would also like to thank Holger P. Kley for useful comments.

\section{Preliminaries}
\label{pre}

A \textit{curve} will always be
reduced and irreducible in this paper.

We now quote some results which will be needed in the rest of the
paper. Most of these results are due to Saint-Donat \cite{S-D}.

\begin{prop}    \label{bp-fix}
  \cite[Cor. 3.2]{S-D} Let $\Sigma$ be a complete linear system on a
  $K3$ surface. Then $\Sigma$ has no base points outside its fixed components.
\end{prop}

\begin{prop} \label{bert}    
  \cite[Prop. 2.6(i)]{S-D} Let $|C| \not = \emptyset$ be a complete linear system without fixed
  components on a $K3$ surface $X$ such that $C^2 > 0$. Then the generic
  member of $|C|$ is smooth and irreducible and $h^1(\O_X(C)) = 0$
\end{prop}

\begin{prop}  \label{ellpen}
  Let $|C| \not = \emptyset$ be a complete linear system without fixed
  components on a $K3$ surface such that $C^2 = 0$. Then every 
  member of $|C|$ can be written as a sum $E_1+E_2+...+E_k$, where $E_i
  \in |E|$ for $i=1,...,k$ and $E$ is a smooth curve of genus $1$.

  In other words, $|C|$ is a multiple $k$ of an elliptic pencil. 
  
  In particular, if $C$ is part of a basis of $\Pic X$, then the generic
  member of $|C|$ is smooth and irreducible. 
\end{prop}

\begin{proof}
  This is \cite[Prop. 2.6(ii)]{S-D}. For the last statement, since $C$
  is part of the basis of $\Pic X$, $k=1$ and $|C| = |E|$. 
\end{proof}

We will also need the following criteria for base point freeness and
very ampleness of a line bundle on a $K3$ surface.

\begin{lemma} \label{vafn}
  \cite{S-D}(see also \cite[Thm. 1.1]{kn1}) Let $L$ be a nef line bundle on a $K3$ surface. Then

(a) $|L|$ is not base point free if and only if there exist curves $E,
    \Gamma$ and an integer $k \geq 2$ such that 
    \[ L \sim kE+ \Gamma, \hspace{.2in} E^2=0, \hspace{.2in}  \Gamma^2=-2,
    \hspace{.2in} E.\Gamma=1. \] 
  In this case, every member of $|L|$ is of the form
  $E_1+...+E_k+\Gamma$, where $E_i \in |E|$ for all $i$. 

  Equivalently, $L$ is not base point free if and only if there is a
  divisor $E$ satisfying $E^2=0$ and $E.L=1$.

(b) $L$ is very ample if and only if $L^2 \geq 4$ and
\begin{itemize}
  \item[(I)]there is no divisor $E$ such that $E^2=0$,
    $E.L=1,2$,
  \item[(II)]there is no divisor $E$ such that $E^2=2$,
    $L \sim 2E$,   and
  \item[(III)]there is no divisor $E$ such that $E^2=-2$,
    $E.L=0$,
  \end{itemize}

\end{lemma}

Note that (II) in (b) is automatically fulfilled if $L$ is a part of a
basis of $\Pic X$, which it often will be in our cases.

Let $L$ be a base point free line bundle on a $K3$ surface with 
$\dim |L|=r \geq 2$. Then $|L|$ defines a morphism
\[ \phi _L : X \hpil \PP^r, \]
whose image $\phi_L (X)$ is called a projective model of $X$. 

We have the following result:

\begin{prop} \label{bircond}
 \cite{S-D} (i) If there is a divisor $E$ such that $E^2=0$ and
    $E.L=2$, or $E^2=2$ and $L \sim 2E$,
    then $\phi_L$ is $2:1$ onto a surface of degree $\frac{1}{2}L^2$.
 
(ii) If there is no such divisor, then $\phi_L$ is birational
    onto a surface of degree $L^2$ (in fact it is an isomorphism
    outside of finitely many contracted smooth rational curves), and 
  $\phi_L (X)$ is normal with only rational double points.

\end{prop}

In \cite{S-D} an $L$ which is base point free and
as in (i) is called \textit{hyperelliptic}, as in this case all
smooth curves in $|L|$ are hyperelliptic. We will call an $L$ which
is base point free and  
as in (ii) \textit{birationally very ample} (see \cite{kn1} for a
generalization). 

We will need the criteria for very ampleness to prove Theorem
\ref{comp.char} and for birational very ampleness to prove the
statement in Remark \ref{mainrem}. 

We also have the following result about the ideal of  $\phi_L (X)$,
when $L$ is birationally very ample.

\begin{prop}  \label{surgen}
  \cite[Thm. 7.2]{S-D} Let $L$ with $H^2 \geq 8$ be a birationally very
  ample divisor on a $K3$ surface $X$. Then the ideal of $\phi_L (X)$
  is generated by its elements of degree $2$, except if there exists a curve
  $E$ such that $E^2=0$ and $E.L=3$, in which case the ideal of
  $X$ is generated by its elements of degree $2$ and $3$. 
\end{prop}

We will concentrate on the proof of Theorem \ref{comp.char}, and give
the main ideas of  the proof of the statement in Remark \ref{mainrem}
in Remark \ref{proofmainrem} below.  

The immediate restrictions on the degree and genus of a divisor come
from the Hodge index theorem. Indeed, if $H$ is any divisor on a $K3$
surface satisfying $H^2=2n >0$ and $C$ is any divisor satisfying $C.H
=d$ and $C^2=2(g-1)$, we get from the Hodge index theorem:
\[ (C.H)^2-C^2H^2 = d^2-4n(g-1) \geq 0, \]
with equality if and only if
\[ dH \sim 2nC. \]

\section{The case $d^2-4n(g-1)=0$}
\label{=0}

We have 

\begin{prop}  \label{dep}
  Let $n \geq 2$, $d>0$, $g \geq 0$ be integers satisfying
  $d^2-4n(g-1) = 0$. Then there is a $K3$ surface $X$ of degree $2n$ in
  $\PP^{n+1}$ containing a smooth curve $C$  of degree $d$ and genus $g$
  if and only if  
  there exist integers $k,m \geq 1$ and $(k,m) \not = (2,1)$ such that
  \[ n = k^2 m  \hspace{.2in} \mbox{and} \hspace{.2in} 2n \hs
  \mbox{divides} \hs kd.  \]
  
  Furthermore, $X$ can be chosen such that $\Pic X = \ZZ
  \frac{2n}{dk} C = \ZZ \frac{1}{k}H$, where $H$ is 
  the hyperplane section of $X$, and if $n \geq 4$, such that $X$ is 
  scheme-theoretically an intersection of quadrics. 
\end{prop}

\begin{proof}
  First we show that these conditions are necessary. Let $X$ be a
  projective $K3$ surface containing a smooth curve of type $(d,g)$ such
  that $d^2-4n(g-1) = 0$. Let $H$ be a hyperplane section of
  $X$. Since $C \sim \frac{d}{2n} H$, there has to exist a divisor
  $D$ and an integer $k \geq 1$ such that $H \sim kD$ and $2n$ divides
  $kd$. Furthermore, letting $D^2=2m$, $m \geq 1$, one gets 
  \[ H^2=2n=2mk^2 ,\]
  so $n = k^2 m$.

  If $(k,m) = (2,1)$, then $n=4$ and $H \sim 2D$ for a divisor $D$ such
  that $D^2=2$, but this is impossible by Lemma \ref{vafn}.

  Now we show that these conditions are sufficient by explicitly
  constructing a projective $K3$ surface of degree $2n$ containing a
  smooth curve of type $(d,g)$ under the above hypotheses.
  
 Consider the  rank $1$ lattice $L=\ZZ D$ with intersection form $D^2 =
 2m$. This lattice is integral and even, and it has signature
 $(1,0)$.  

 Now \cite[Thm. 2.9(i)]{morr} (see also \cite{niku}) states that: 

 \textit{If $\rho \leq 10$, then every even lattice of signature $(1,
  \rho -1)$ occurs as the N\'{e}ron-Severi group of some algebraic $K3$
  surface.}

 Therefore there exists an algebraic $K3$ surface $X$ such that $\Pic X =
 \ZZ D$. Since $D^2>0$, either $|D|$ or $|-D|$ contains an effective
 member, so 
 we can assume $D$ is effective (possibly after having changed $D$ with
 $-D$). Since $D$ generates $\Pic X$, $D$ is ample by Nakai's criterion
 (in particular, it is nef). If $m \geq 2$, then $D^2 \geq 4$, and by
 Lemma \ref{vafn} $D$ is very ample. Indeed, since $\Pic X = \ZZ D$ and
 $D^2 >0$, there can exist no divisor $E$ such that $E^2 = 0, -2$ or $D
 \sim 2E$.

 If $m=1$, then by our assumptions $k \geq 3$, so by 
 \cite[Thm. 8.3]{S-D} $kD$ is very ample. 

 So in all cases, $H := kD$ is very ample under the above assumptions,
 and by $H^2 = k^2 D^2= 2mk^2 = 2n$, we can embed $X$ as a $K3$ surface
 of degree $2n$ in $\PP^{n+1}$. Define
 \[ C := \frac{dk}{2n} D. \]
 then $C$ is nef, since it is a non-negative multiple of a nef divisor, and
 by Lemma \ref{vafn} it is base point free (no curves with
 self-intersection $-2$ or $0$ can occur on $X$ since $D$ generates the
 Picard group).
 
 So by Propositions \ref{bp-fix} and \ref{bert} the generic member of
 $|C|$ is smooth and irreducible, and one easily checks that $C.H=d$ and
 $C^2=2g-2$.

 The two last assertions follow from the construction (the last one
 follows from Proposition \ref{surgen}).
\end{proof}

\section{The case $d^2-4n(g-1)>0$}
\label{>0}

We have the following result of Oguiso:

\begin{thm} \label{og2}
  Let $n \geq 2$ and $d \geq 1$ be positive integers. Then there exist
  a $K3$ surface $X$ of degree $2n$ in $\PP ^{n+1}$ and a smooth
  rational curve $C$ of degree $d$ on $X$.

  Furthermore, $X$ can be chosen such that 
  $ \Pic X = \ZZ H \+ \ZZ C$, where $H$ is  
  the hyperplane section of $X$, and if $n \geq 4$, then $X$ can be
  chosen to be 
  scheme-theoretically an intersection of quadrics. 
\end{thm}

\begin{proof}
  This is \cite[Thm 3]{oguiso}. The last statement follows again by
  Proposition \ref{surgen}  since $| \disc (H,C) | = d^2+4n >16$ and
  a divisor $E$ such as the one in the proposition would give
  $|\disc (E,H)| = 9 $.
\end{proof}

We can make a more general construction:

\begin{prop}   \label{const}
  Let $n \geq 1$, $d \geq 1$, $g \geq 0$ be positive integers satisfying
  $d^2-4n(g-1) > 0$. Then there exist an algebraic $K3$ surface $X$ and two
  divisors $H$ and $C$ on $X$ such that $\Pic X = \ZZ H \+ \ZZ C$,
  $H^2=2n$, $C.H = d$, $C^2=2(g-1)$ and $H$ is nef.
\end{prop}

\begin{proof}
  Consider the lattice $L=\ZZ H \+ \ZZ C$ with intersection matrix 
\[  \left[ 
  \begin{array}{cc}
  H^2  &  H.C   \\ 
  C.H  &  C^2 
    \end{array} \right]  = 
    \left[
  \begin{array}{cc}
  2n  &  d   \\ 
  d   &  2(g-1) 
    \end{array} \right]     \]  
This lattice is integral and even,
and it has signature $(1,1)$ if and only if $d^2-4n(g-1) > 0$. 

By \cite[Thm. 2.9(i)]{morr} again, we conclude that the lattice $L$
occurs as the Picard group of some algebraic $K3$ surface $X$. It remains to
show that $H$ can be chosen nef.  

Consider the group generated by the Picard-Lefschetz reflections
\begin{eqnarray*}
      \Pic X  & \stackrel{\pi_{\Gamma}}{\hpil} &  \Pic X  \\
        D     & \longmapsto           & D+(D.\Gamma)\Gamma,
\end{eqnarray*}  
where $\Gamma \in  \Pic X$ satisfies $\Gamma ^2 = -2$ and $D \in  \Pic
X$ satisfies $D^2>0$. Now \cite[VIII, Prop. 3.9]{BPV} states that
a fundamental domain for this action is the \textit{big-and-nef cone} of
$X$. Since $H^2>0$, we can assume that $H$ is nef.
\end{proof}

We would like to investigate under which conditions $H$ is very ample and
$|C|$ contains a smooth irreducible member. To show the latter for
$g>0$, it will be enough to show that $|C|$ is base point free, by
Propositions \ref{bert} and \ref{ellpen}.
 
We first need a basic lemma.

\begin{lemma} \label{basiclemma}
  Let $H$, $C$, $X$, $n$, $d$ and $g$ be as in Proposition
  \ref{const} and $k \geq 1$ an integer. If $(d,g)=(2nk,nk^2)$ (resp.
  $(d,g)=(nk, \frac{nk^2+3}{4})$), 
  we can assume (after a change of basis of $\Pic X$) that $kH-C >0$
  (resp. $kH-2C >0$).
\end{lemma}

\begin{proof}
  We calculate $(kH-C)^2=-2$ (resp. $(kH-2C)^2=-2$), so by Riemann-Roch
  either $kH-C >0$ or $C-kH >0$ (resp. either $kH-2C >0$ or $2C-kH >0$).
  If the latter is the case, define $C':=2kH-C$ (resp. $C':=kH-C$). Then
  one calculates $C'.H=d$ and ${C'}^2=2(g-1)$, and since clearly 
  $\Pic X \iso \ZZ H \+ \ZZ C'$, we can substitute $C$ with $C'$.
\end{proof}

\begin{prop}  \label{amim}
  Let $H$, $C$, $X$, $n$, $d$ and $g$ be as in Proposition \ref{const}
  with $g \geq 1$, and with the additional assumptions
  that $kH-C >0$ (resp. $kH-2C >0$) if $(d,g)=(2nk,nk^2)$ (resp. if 
  $(d,g)=(nk,\frac{nk^2+3}{4})$).
  Assume $H$ is base point free. Then $|C|$ contains a smooth irreducible
  member if and only if we are not in one of the following cases:
  \begin{itemize}
  \item[(i)]  $(d,g) = (2n+1,n+1)$, 
  \item[(ii)] $d^2-4n(g-1)=1$ and $d-1$ or $d+1$ divides $2n$.
  \end{itemize}
\end{prop}

\begin{proof}
  We first show that $C$ is nef except for the case (i).

  Assume that $C$ is not nef. Then there is a
  curve $\Gamma$ (necessarily contained in the fixed component of $|C|$)
  such that $C.\Gamma <0$ and $\Gamma^2=-2$. 
 
  We now consider the two cases $\Gamma.H >0$ and $\Gamma.H=0$
  \footnote{This latter case occurs only if $H$ is not ample, so it is
  only interesting in order to prove the statement in Remark
  \ref{mainrem}. To prove Theorem \ref{comp.char} we could assume that
  $H$ is ample and thus get an easier proof of Proposition
  \ref{amim}.}.

  If $\Gamma.H >0$, define $a:= - C.\Gamma \geq 1$ and 
  \[  C' := C -a \Gamma  ,\]
  then 
  \[ {C'}^2 = 2(g-1) \geq 0,\]
  so by Riemann-Roch either $|C'|$ or $|-C'|$ contains an effective
  member, and since $h^0(\O_X (a \Gamma)) =1$, it must be $|C'|$. Hence
\[  0 < d':= C'.H = C.H -a (\Gamma.H)  < d , \]
  where we have used that $\Gamma.H>0$ to get the strict inequality on
  the right and ${C'}^2 \geq 0$ to get the strict inequality on the
  left, by the Hodge index theorem. Since $H^2 >0$, one 
  must have ${d'}^2-4n(g-1) \geq 0$ by the Hodge index theorem, and
  equality occurs if and only if $d'H \sim 2nC'$. 
  
  We now show that ${d'}^2-4n(g-1) = 0$ only if $(d,g) = (2n+1,n+1)$
  and that $C$ is not nef in this case.

  Write $\Gamma \sim xH + yC$, for two integers $x$ and $y$. We have
\[ C' \sim C - a\Gamma \sim C - a(xH+yC) \sim -axH + (1-ay)C \sim
\frac{d'}{2n}H, \]
  which implies $y=a=1$. We then have
\[ -1 = \Gamma.C = dx + 2(g-1)y = dx + 2(g-1), \] 
  which yields $x= - \frac{2g-1}{d}$, whence
\[ \Gamma \sim - \frac{2g-1}{d} H + C. \]
  Note that this implies 
  \begin{equation}
    \label{eq:first}
    d \hs | \hs 2g-1.
  \end{equation}
We now use 
\[ -2 = \Gamma ^2 = \frac{(2g-1)^2}{d^2}2n-2(2g-1)+2(g-1) = 2
(\frac{(2g-1)^2 n}{d^2}-g) \]
to conclude
\begin{equation}
  \label{eq:second}
  n= \frac{(g-1) d^2}{(2g-1)^2}.
\end{equation}
Using this, we calculate
\[ \Gamma.H = - 2n \frac{2g-1}{d} +d = - 2 \frac{(g-1) d^2}{(2g-1)^2}
\frac{2g-1}{d}+d = \frac{d}{2g-1} , \]
which yields
\begin{equation}
    \label{eq:third}
    2g-1\hs | \hs d.
  \end{equation}
Comparing (\ref{eq:first}) and (\ref{eq:third}), we get $d=2g-1$,
which gives by (\ref{eq:second}) that
\[ (d,g) = (2n+1,n+1). \]

So we have shown that ${d'}^2-4n(g-1) = 0$ occurs only when $(d,g) =
  (2n+1,n+1)$ and that $C$ is not nef in this case.

We now consider the case when ${d'}^2-4n(g-1) > 0$. Since 
  \[ 0 \not = {d'}^2-4n(g-1) = |\disc(H,C')| < d^2-4n(g-1)=
      |\disc(H,C)|,\] 
  then $\disc(H,C)$ cannot divide $\disc(H,C')$ and we have a
  contradiction, so $C$ is nef.

  If $\Gamma.H =0$, write $\Gamma \sim xH+yC$. We have 
\[ -2 = \Gamma ^2 = \Gamma.(xH+yC) = yC.\Gamma, \]
  which gives the two possibilities
  \begin{itemize}
  \item[(a)] $y=1$, $C.\Gamma=-2$, and
  \item[(b)] $y=2$, $C.\Gamma=-1$.
  \end{itemize}

In case (a), we get from $\Gamma.H =2nx+dy = 2nx+d=0$ that $x =-d/2n$,
which means that
\[ d=2nk \hs \mbox{ and } x=-k, \]
for some integer $k \geq 1$.
From $\Gamma.C = dx + 2(g-1)y = -2nk^2 +2(g-1)=-2$, we get $g=nk^2$. So
$(d,g)=(2nk,nk^2)$ and $\Gamma \sim -kH +C$, which by 
assumption is not effective, a contradiction.

In case (b), we get from $\Gamma.H =2nx+dy = 2nx+2d=0$ that $x =-d/n$,
which means that
\[ d=nk \hs \mbox{ and } \hs x=-k, \]
for some integer $k \geq 1$.
We get from $\Gamma.C = dx + 2(g-1)y = -nk^2 + 4(g-1)=-1$, that $g =
(nk^2+3)/4$. So $(d,g)=(nk, \frac{nk^2+3}{4})$ and $\Gamma \sim -kH +2C$,
which by assumption is not effective, again a contradiction. 

So we have proved that $C$ is nef except for the case (i).

If $|C|$ is not  base point free, then using Lemma \ref{vafn}(a),
  $X$ must contain two divisors $E$ and $\Gamma$ such that $E^2=0$ and 
  \[ | \disc (E,\Gamma)| = 1 , \]
  and this must be divisible by $|\disc(H,C)| = d^2-4n(g-1)$, which then
  must be equal to $1$.

  Setting $E \sim xH+yC$ one finds from $E.C=dx+2(g-1)y=1$ and $E^2=2n
  x^2+ 2dxy + 2(g-1)y^2=0$ that
  \[ x = \pm 1 \hspace{.2in}  \mbox{and} \hspace{.2in} y = \frac{1 \mp
    d}{2(g-1)}.\] 
  Using the fact
  that $d^2-4n(g-1)=1$, we get $2(g-1)=\frac{(d+1)(d-1)}{2n}$, which
  gives
\[ (x,y)= (1, - \frac{2n}{d+1}) \hs \mbox{ or } \hs (-1, \frac{2n}{d-1}). \] 

  So if $d^2-4n(g-1)=1$ and $d+1$ or $d-1$ divides $2n$, then 
the divisor $E$ will satisfy $E^2=0$ and $E.C=1$, whence $C$ is not
base point free by Lemma \ref{vafn}(a).

This concludes the proof of the proposition.  
\end{proof}

Note that we have also proved 

\begin{cor}
  Let $H$ and $C$ be divisors on a $K3$ surface $X$ such that $H$ is nef,
  $H^2=2n$, $C.H=d$ and $C^2=2(g-1)$ for some integers $n \geq 1$, $d
  >0$ and $g \geq 1$. If either
  \begin{itemize}
  \item[(a)] $(d,g) = (2n+1,n+1)$, or
  \item[(b)] $d^2-4n(g-1)=1$ and $d+1$ or $d-1$ divides $2n$,  
   \end{itemize}
  then $C$ is not base point free.
 
  In particular, a projective $K3$ surface of degree $2n$, for an
  integer $n \geq 2$ (or even a birational
  projective model of a $K3$ surface) cannot contain an effective,
  irreducible divisor of
  degree $d$ and arithmetic genus $g$ for any values of $d$ and $g$ as
  in (a) or (b). 
\end{cor}

\begin{proof}
  In case (a) the divisor $\Gamma:= C-H$ is effective and satisfies
  $\Gamma.C=-1$, so $C$ is not even nef.

  In case (b) the divisor $E:= H-\frac{2n}{d+1}C$ or
  $E:=-H+\frac{2n}{d-1}C$ is effective and satisfies $E^2=0$ and $E.C=1$.
  Thus $C$ is not base point free by Lemma \ref{vafn}(a).
  \end{proof}

One gets the following 

\begin{thm}  \label{fullcon}
  Let $n \geq 2$, $d \geq 1$, $g \geq 0$ be positive integers satisfying
  $d^2-4ng > 0$ and $(d,g) \not = (2n+1,n+1)$. Then there exists a
  projective $K3$ surface $X$ of degree 
  $2n$ in $\PP^{n+1}$ containing a smooth curve $C$ of degree $d$ and
  genus $g$. Furthermore, we can find an $X$ such that $\Pic X = \ZZ H \+ \ZZ C$, where $H$ is the
  hyperplane section of $X$, and $X$ is scheme-theoretically an
  intersection of quadrics  if $n \geq 4$. 
\end{thm}

\begin{proof}
  The case $g=0$ is Theorem \ref{og2}, so we can assume $g>0$. Since
  $|\disc(H,C)| = d^2-4n(g-1) > 4n$, the $H$ constructed as in 
  Proposition \ref{const} is very ample, since the existence of such
  divisors as in (I) and (III) in Lemma \ref{vafn}(b)
  implies that $| \disc (H,C) |$ must divide $| \disc (H,E) | =1,4,4n$,
  respectively. Now Proposition \ref{amim} gives the rest. 

  To prove that $X$ is an
  intersection of quadrics when $n \geq 
  4$, by Prop \ref{surgen} it is sufficient to show that there
  cannot exist any divisor $E$ such that $E^2=0$ and $E.H=3$. 

  Such an $E$ would give
  \[  |\disc (E,H)| = 9 \]  \label{genquad}  
  but 
  we have $|\disc(H,C)| >4n \geq 16$, when $n \geq 4$.
\end{proof}

Now we only have to investigate the pairs $(d,g)$ for which 
\[ 0< d^2-4n(g-1) \leq 4n  .\]

We proceed as follows. For given $n,d,g$ we use the
construction of Proposition \ref{const} and then investigate whether $H$ is
very ample by using Lemma \ref{vafn}. Then two cases may occur:

\begin{enumerate}
\item Using the fact that $\Pic X = \ZZ H + \ZZ C$ and $H^2=2n$,
  $C.H=d$, $C^2=2(g-1)$, we find that there cannot exist any divisor
  $E \sim xH+yC$ 
  as in cases (I) and (III) of Lemma \ref{vafn}(b), so $H$
  is very ample and by Proposition \ref{amim}, $|C|$ contains a smooth
  irreducible member and there exists a projective $K3$ surface $X$ of degree 
  $2n$ in $\PP^{n+1}$ containing a smooth curve $C$ of degree $d$ and 
  genus $g$.
\item Using the numerical properties $H^2=2n$,
  $C.H=d$ and $C^2=2(g-1)$, we find  a divisor $E \sim aH+bC$ for $a,b
  \in \ZZ$ satisfying case (I) or
  (III) of Lemma \ref{vafn}(b), thus contradicting the very ampleness
  of $H$. This then 
  implies that there cannot exist any projective $K3$ surface of
  degree $2n$ in $\PP^{n+1}$ containing a divisor of degree $d$ and
  genus $g$.
\end{enumerate}

(To prove the statement in Remark \ref{mainrem}, we proceed in an
analogous way, but check the conditions for $H$ to be birationally
very ample instead. We then have to investigate whether any of the
smooth curves in $|C|$ are mapped isomorphically to a smooth curve.
See Remark \ref{proofmainrem} below.) 

For any triple of integers $(n,d,g)$ such that $n \geq 2$, $d >0$,
$g \geq 1$ and $0 < d^2-4n(g-1) \leq 4n$, define 
\[ \k (n,d,g) := d^2-4n(g-1) = | \disc (H,C) |   ,\]
We check conditions (I) and
(III) in Lemma \ref{vafn}. We let 
\[ E \sim xH+yC  ,\]
and use the values of $E^2=2nx^2 + 2dxy + 2(g-1)y^2$ and $E.H = 2nx +dy$
to find the integers $x$ and $y$ (if any).

We get the two equations
\[ n \k (n,d,g)x^2 - (E.H) \k (n,d,g) x - (g-1) (E.H)^2 + \frac{d^2}{2} E^2 =0 \]
and
\[ y = \frac{(E.H) - 2nx}{d}.\]

\begin{itemize}
\item[(a)] If $E.H =1$ and $E^2 = 0$, then 
  \[ x= \frac{\k (n,d,g) \pm d \sqrt{\k (n,d,g)}}{2n \k (n,d,g)},
  \hspace{.5in} y = \mp 
  \frac{1}{\sqrt{\k (n,d,g)}}, \]
  and the only possibility is  $\k (n,d,g)=1$, so 
  \[ x= \frac{1 \pm d}{2n},\hspace{.5in} y = \mp 1,\]
  and we must have $d \eqv \pm 1 \hs (\mod 2n)$.
\item[(b)] If $E.H =2$ and $E^2 = 0$, then 
  \[ x= \frac{\k (n,d,g) \pm d \sqrt{\k (n,d,g)}}{n\k}, \hspace{.5in} y = \mp
  \frac{2}{\sqrt{\k (n,d,g)}}, \]
  and the only possibilities are  $\k (n,d,g)=1$ or $4$.

  If  $\k (n,d,g)=1$, then 
  \[ x= \frac{1 \pm d}{n},\hspace{.5in} y = \mp 2,\]   
  and we must have $d \eqv \pm 1 \hs (\mod 2n)$ or 
  $d \eqv n \pm 1 \hs (\mod 2n) $. 

  If  $\k (n,d,g)=4$, then 
  \[ x= \frac{2 \pm d}{2n},\hspace{.5in} y = \mp 1,\]   
  and we must have $d \eqv \pm 2 \hs (\mod 2n)$.

\item[(c)] Finally, if $E.H =0$ and $E^2 = -2$, we get 
\[ -2 =E^2 = E.(xH+yC) = yE.C, \]
so $y=-1$ or $-2$ (since $C$ is nef), and by $E.H=2nx+dy=0$, we get 
\[ d \eqv 0 \hs (\mod 2n) \hs \mbox{ and } \hs x= \frac{d}{2n}, \hs \hs 
    \mbox{ or } \hs \hs d \eqv 0 \hs (\mod n) \hs \mbox{ and } \hs x=
    \frac{d}{n} \] 
respectively.

One now easily calculates (using $E.C=2$ (resp. $1$))
\[ \k (n,d,g)=4n \hs \mbox{ and } \hs n \]
respectively. Furthermore, in the latter case, if $d \eqv 0 \hs (\mod
2n)$, writing $d=2nk$, for some integer $k \geq 1$, we get from
$E.C=1$ the absurdity $g =nk^2+ \frac{3}{4}$. So
we actually have $d \eqv n \hs (\mod 2n)$ in this case.
\end{itemize}

What we have left to prove in Theorem \ref{comp.char} is that $X$ can be
chosen as an intersection of quadrics in case (iii) and under the given
assumptions in case (ii).

We have to show, by Proposition \ref{surgen}, that there cannot exist
any divisor $E$ such that $E.H=3$ and $E^2=0$. Since such an $E$ would
give $| \disc (E,H)| = 9$, and we have $| \disc (H,C)| = d^2-4n(g-1) =
4n \geq 16$ if $n \geq 4$ in case (iii), this case is proved.

For case (ii), we proceed as above and set $E \sim xH+yC$, and try to
find the integers $x$ and $y$. We find
\[ x= \frac{3 (\k (n,d,g) \pm d \sqrt{\k (n,d,g)})}{2n \k (n,d,g)}, \hspace{.5in} y = \mp
\frac{3}{\sqrt{\k (n,d,g)}}, \]
so we must have $\k (n,d,g)=1$ or $9$.

If $\k (n,d,g)=1$, then 
\[ x= \frac{3 (1 \pm d)}{2n}, \hspace{.5in} y = \mp 3, \]
and $E$ exists if and only if $3d \eqv  \pm 3 \hs (\mod 2n)$.

If $\k (n,d,g)=9$, then 
\[ x= \frac{3 \pm d}{2n}, \hspace{.5in} y = \mp 1, \]
and $E$ exists if and only if $d \eqv  \pm 3 \hs (\mod 2n)$.

So, under the given assumptions, the constructed $X$ is an intersection
of quadrics. Conversely, given a $K3$ surface $X$ of degree $2n$
containing a smooth curve $C$  of degree $d$ and genus $g$ such that
$d^2-4(n-1)=1$ and $3d \eqv  \pm 3 \hs (\mod 2n)$ or  
$d^2-4(n-1)=9$ and $d \eqv  \pm 3 \hs (\mod 2n)$, then the divisor $E$
above, which is a linear combination of $C$ and the hyperplane section
$H$, will be a divisor such that $E.H=3$ and $E^2=0$. Hence $X$ must be
an intersection of both quadrics and cubics. 

This concludes the proof of Theorem \ref{comp.char}.

\begin{rem} \label{proofmainrem}
  {\rm If we allow $\phi_H (X)$ to be a birational projective model of
  $X$ (i.e. we require $H$ to be birationally very ample only), we can
  allow the existence of a divisor $E$ such that $E^2=-2$ and $E.H=0$.
  That is, we can allow the cases
\begin{itemize}
\item[] $d^2-4n(g-1)=n$  and $d \eqv n \hs (\mod 2n)$,  \hspace{1cm} and
\item[] $d^2-4n(g-1)=4n$ and $d \eqv 0 \hs (\mod 2n)$ 
\end{itemize}
In the first case, with a lattice as in Proposition  \ref{const},
  the only contracted curve $\Gamma$ satisfies $\Gamma.C=1$, whence
  every smooth curve in $|C|$ is mapped isomorphically by $\phi_H$ to
  a smooth curve of degree $d$ and genus $g$.

  In the second case, define
\[ \Gamma := \frac{d}{2n}H-C. \]
  Then $\Gamma ^2=-2$, $\Gamma.H=0$ and $\Gamma.C=2$, so any
  irreducible member of $|C|$ contains a length two scheme where $H$
  fails to be very ample. Thus $\phi_H (C)$ is singular for all
  irreducible $C'\in |C|$.

  This shows the statement in Remark \ref{mainrem}. } 
\end{rem}

\section{Proof of Proposition \ref{extra} }

This section is devoted to the proof  of Proposition \ref{extra}. 

Let $C'\in |C|$.
By the exact sequence 
\[ 0 \khpil \O_X (kH-C) \khpil \O_X (kH) \khpil \O_{C'} (kH) \khpil 0, \]
and using that $H^1 (\O_X (kH)) = H^2 (\O_X (kH))=0$
by the Kodaira vanishing theorem, we see that (using Serre duality)
\[ h^1 (\O_{C'} (k)) = h^0 (C-kH). \]

If $d \leq 2nk$, then $(C-kH).H =d-2nk \leq 0$, which implies $h^0
(C-kH)=0$, since $H$ is ample.

If $d > 2nk$ and $dk \leq nk^2+g$, we have $(C-kH)^2 \geq -2$ and
$(C-kH).H >0$. So by Riemann-Roch, $C-kH >0$.

To finish the proof, let $d > 2nk$ and $dk > nk^2+g$ and assume that
there is an element $D \in |C-kH|$. Then $D^2 < -2$, so
$D$ has to contain an irreducible curve $\Gamma$ such that $D.\Gamma
<0$ and $\Gamma^2=-2$. As seen above, we can assume that $\Pic X$ is
either generated by some rational multiple of the hyperplane section
or generated by $H$ and $C$. Clearly, in the first case, all divisors
have non-negative self-intersection, so we can assume 
$\Pic X = \ZZ H \+ \ZZ C$. 

We can write 
\[  D = m \Gamma + E,  \]
with $\Gamma$ not appearing as a component of $E$ and $E \geq 0$. 

If $E=0$, then $m=1$ since $D$ is a part of a basis of $\Pic X$, but
then $D^2= -2$, which is a contradiction, so $E > 0$. We then see that
$D.\Gamma \geq -2m$ and $\Gamma.H < D.H/m$.

Now we define
the divisor \[  D' := D + (D.\Gamma)\Gamma.\]
Then ${D'}^2=D^2$ and 
\[ -D.H < D'.H = D.H + (D.\Gamma)(\Gamma.H) < D.H. \]
By the Hodge index theorem and the fact that ${D'}^2 <0$, we get ${D'}^2 H^2 < (D'.H)^2$, so we have
\[ 0 \not = (D'.H)^2- {D'}^2 H^2 = | \disc (H,D')| <  (D.H)^2- D^2 H^2 =
| \disc (H,D)|, \]
a contradiction, since clearly $\Pic X = \ZZ H \+ \ZZ D$.  

This concludes the proof of Proposition \ref{extra}.

\section{Application to complete intersection $K3$ surfaces}
\label{app}
 
$K3$ surfaces of degree $4$ in $\PP^3$ or of degree $6$ in $\PP^4$ have to
be smooth quartics and smooth complete intersections of type $(2,3)$
respectively. Furthermore, by Proposition \ref{surgen} a complete intersection $K3$ surface of
degree $8$ in $\PP^5$ has to be either an intersection of quadrics, in
which case it is easily seen to be a complete intersection of type
$(2,2,2)$, or an intersection of both quadrics and cubics, in which case
it cannot be a complete intersection.

Applying Theorem  \ref{comp.char} to the three kinds of complete intersection
$K3$ surfaces, one gets:

\begin{thm} \label{genmori}
  Let $d>0$ and $g \geq 0$ be integers. Then:
  \begin{enumerate}
  \item (\textbf{Mori \cite{mori} })There exists a smooth quartic
    surface $X$ in 
    $\PP^3$ containing a smooth curve $C$ of degree $d$ and genus $g$
    if and only if  a) $g=d^2/8+1$, or  b) $g<d^2/8$ and $(d,g)
    \not = (5,3)$. Furthermore, a) holds if and only if $\O_X (1)$ and
    $\O_X (C)$ are dependent in $\Pic X$, in which case $C$ is a
    complete intersection of $X$ and a hypersurface of degree $d/4$.
  \item There exists a $K3$ surface $X$ of type $(2,3)$ in  $\PP^4$
    containing a smooth curve $C$ of degree $d$ and genus $g$ if and
    only if  a) $g=d^2/12+1$, b) $g=d^2/12+1/4$ or  c) $g<d^2/12$ and $(d,g )
    \not = (7,4)$. Furthermore, a) holds if and only if $\O_X (1)$ and
    $\O_X (C)$ are dependent in $\Pic X$, in which case $C$ is a
    complete intersection of $X$ and a hypersurface of degree $d/6$.
  \item There exists a $K3$ surface $X$ of type $(2,2,2)$ in  $\PP^5$
    containing a smooth curve $C$ of degree $d$ and genus $g$  if and
    only if  a) $g=d^2/16+1$ and $d$ is divisible by $8$,  b) $g=d^2/16$
    and $d \eqv 4  \hs (\mod 8)$, or  c) $g<d^2/16$ and $(d,g )
    \not = (9,5)$. Furthermore, a) holds if and only if $\O_X (1)$ and
    $\O_X (C)$ are dependent in $\Pic X$, in which case $C$ is a
    complete intersection of $X$ and a hypersurface of degree $d/8$.
  \end{enumerate}
\end{thm}

\vspace{.5cm}

\noindent \textit{Mathematics Subject Classification:} 14J28 (14H45).


\noindent \textit{E-mail:} andreask@mi.uib.no

\end{document}